\documentclass{article}

\usepackage[left=2cm,right=1.5cm,top=1.8cm,bottom=1.8cm]{geometry}
\usepackage{multicol}
\usepackage[hidelinks]{hyperref}
\usepackage[T1]{fontenc}

\usepackage[UKenglish]{babel}
\usepackage{csquotes}

\usepackage[labelfont=bf]{caption}
\usepackage{graphicx}
\newenvironment{myfigure}
  {\par\medskip\noindent\minipage{\linewidth}}
  {\endminipage\par\medskip}

\usepackage{titlesec}
\titleformat{\section}
  {\titlerule\centering\normalfont\large\bfseries}{}{1em}{}[\hrule]
\titlespacing{\paragraph}{%
  0pt}{
  0.5\baselineskip}{
  1em}

\usepackage{bibentry}
\newcommand{\MR}{MathRepo }

\usepackage{xpunctuate}
\providecommand{\eg}[0]{e.g\xperiod}

\usepackage{xcolor}

\definecolor{burgundy}{rgb}{0.5, 0.0, 0.13}
\newcommand{\MYhref}[3][]{\href{#2}{\color{#1}{#3}}}%

\author{Claudia Fevola and Christiane G\"orgen\\
\texttt{fevola@mis.mpg.de} and \texttt{goergen@math.uni-leipzig.de}}
\title{\bfseries The mathematical research-data repository \MR}
\date{}

\begin{document}
\maketitle
\begin{multicols}{2}

\section{Research data in mathematics}

MathRepo, located at \texttt{\MYhref[burgundy]{https://mathrepo.mis.mpg.de}{https://mathrepo.mis.mpg.de}}, is an online repository for mathematical research data. Research data, broadly speaking, can be defined as \enquote{the recorded factual material commonly accepted in the scientific community as necessary to validate research findings}\footnote{OMB Circular 110, \url{https://www.whitehouse.gov/sites/whitehouse.gov/files/omb/circulars/A110/2cfr215-0.pdf}}. 
In mathematics, research data comes in many different flavours. For instance, in computer algebra it most often takes the form of mathematical documents, notebooks, research-software packages and libraries, computer algebra systems, algorithms, and collections of mathematical objects \cite[e.g.]{DMVM}. \MR contains foremostly computer-algebra research data of three main types: computations performed for paper publications, additional lists of examples to theoretical results, and presentations of problems solved in workshops or lectures. It provides one central location for mathematicians to store and share any additional data they might want to make publicly available alongside their paper publications.

The repository was established in 2017 at the Max Planck Institute for Mathematics in the Sciences\footnote{\url{https://www.mis.mpg.de}} (MPI MiS) in Leipzig, following the initiative of Bernd Sturmfels, \mbox{Ronald} Kriemann, and Yue Ren. It has had seven different maintainers over the past five years: Yue Ren and Mahsa Sayyary in 2017 and 2018, Lukas K\"uhne and Verena Morys from 2019, and the authors of this paper together with Carlos Am\'endola since 2021.
 There are of course a variety of other storage solutions across different scientific fields\footnote{See the list of research-data repositories at \url{https://www.re3data.org}, \eg}, some research data can be published as software packages \cite[e.g.]{smallgroups}, and some mathematical libraries have their own homepages\footnote{See the library \enquote{Small Phylogenetic Trees} {\url{https://www.coloradocollege.edu/aapps/ldg/small-trees/small-trees_0.html}} or \enquote{The Markov Bases Database} \url{https://markov-bases.de}, \eg}. However, an infrastructure which can both store and visually present a wide range of different data in a wiki- or blog-entry style was---and at the time of writing \emph{is}---not yet broadly established. In particular, for smaller contributions to papers like proofs by computation in a particular programming language or lists of examples of objects with properties of interest, a centralised infrastructure was completely missing. \MR aimed to fill this gap.

At the time of writing, the repository has gathered a total of forty individual contributions. The contributor community is the nonlinear-algebra working group\footnote{\url{https://www.mis.mpg.de/nlalg/nlalg-people.html}} at MPI MiS, though being a member is not necessary to gain access to the repository. Informally registering with MPI MiS's IT service is sufficient to obtain contributer rights. \MR is hosted on the servers of MPI MiS\footnote{\url{https://gitlab.mis.mpg.de}}. It is planned to last for at least the coming decade and is limited in size and capacities as is the underlying open-source software \texttt{GitLab}\footnote{\url{https://docs.gitlab.com/ee/administration/instance_limits.html}}. Just like with paper-storage options such as \url{arXiv.org}, also in \MR there is no hard quality control. Merge requests are accepted after a brief visual check of three basic requirements: that the content of a new page is mathematical research data, that it provides references to relevant literature, and that the authors are generally known in the nearby scientific community. In particular, there are no strict rules for the actual presentation of the mathematical content. Using \MR as a pure repository is just as possible as using it to provide an in-depth introduction to a topic of interest, as we will see in a number of examples below. 

\section{What makes a good repository?}

In 2021 the Mathematical Research Data Initiative MaRDI\footnote{\url{https://www.mardi4nfdi.de}} has set about developing and establishing infrastructure for research data in the German mathematics community. 
The declared purpose of the consortium is to establish the FAIR principles for these data: their long-term findability, accessability, interoperability, and reusability shall be ensured \cite{FAIR,NFDI}. Currently, twenty-five partnering research organisations within MaRDI are working towards this aim. These are universities and institutes from the Fraunhofer and Max Planck societies as well as the Leibniz Gemeinschaft, the professional societies DMV, GAMM, and GOR, the European Mathematical Society and partners from mathematically-inclined Clusters of Excellence within Germany. This network ensures a nationally consistent implementation and follows a bottom-up approach in setting up new standards. In MaRDI's computer-algebra task area, the \texttt{OSCAR}\footnote{\url{https://oscar.computeralgebra.de}} group is a key player, for instance.

One way of implementing the FAIR principles in practice is via the usage of trustworthy storage solutions for research data. This article is a snapshot of the service that \MR provides to the computer-algebra and other mathematical communities in this context and at this precise point in time, in February 2022---after nearly five years of maintenance under different leads, with new programming languages and software solutions coming up, and MaRDI waiting in the wings.

All present and past maintainers of \MR have themselves contributed to the repository. In particular, the initiators had plenty hands-on experience in handling \emph{not} FAIR research data: data which was promised in papers and stored on long-gone personal homepages, data blocked by pay walls, data which would run on one computer but not on another, and data which would just not provide the promised results. Their key strategy to address these issues with \MR was usability. They envisioned that if the repository was easy to access for authors that would offer a practical counterpart to other, decentralised, research-type specific or hard-to-maintain storage options. And if it was used by authors, it would become known to readers as well and thus self-establish in the community, replacing cumbersome past solutions. A low entry barrier was initially achieved by Yue Ren being the sole maintainer who uploaded and curated all of the contributions in 2017 and '18. It was later replaced by annual \texttt{GitLab}-training sessions with local contributors at the MPI MiS. 
Both strategies have had success and facilitated acceptance of the service in their local academic community, though \MR has never been aimed solely at that particular audience. 

The key question for us is now: is the content of \MR FAIR in a MaRDI context? Or, more practically for you as the reader, is using \MR for your own research a solution for the future? We will discuss these questions over the coming sections, presenting the mathematics currently present on the repository and discussing what the FAIR principles mean for these in practice.

\section{The mathematics in \MR}
The common theme of all research data currently present on \MR is \emph{nonlinear algebra}. This is a diverse and developing field of mathematics promoted by the recent expansion of nonlinear methods across applications \cite{michalek2021invitation,sturmfels2021beyond}. Of course, the theory, algorithms, and software from linear algebra and numerical linear algebra have a crucial function in the process of modelling problems arising in the natural sciences and engineering. But the natural occurrence of nonlinear equations in real-world applications together with an increasingly strong toolbox of new computational methods brought about a growing use of nonlinear approaches to mathematical modelling. These recent developments rely on techniques from algebraic geometry, topology, combinatorics, group theory, commutative algebra or representation theory.
Vice versa, applications are also a central motor for driving new results in this field. Examples of this can be found in physics, polynomial optimization, partial differential equations, algebraic statistics, and algebraic vision \cite{breiding2021nonlinear}.

In practice, many of the applications of nonlinear algebra boil down to the problem of finding solutions to systems of polynomial equations. Broadly speaking, there are then two main computational approaches: symbolic and numerical. The former often relies on Gröbner bases' computations, the latter typically uses homotopy continuation. 
Several software and programming languages provide effective tools for these computations and \MR furnishes a valid collection of the most well known ones, including \texttt{GAP}, \texttt{Julia}, \texttt{Macaulay2}, \texttt{Magma}, \texttt{Maple}, \texttt{Mathematica}, \texttt{Matlab}, \texttt{Polymake}, 
\texttt{Sage}, and \texttt{Singular}. 

\paragraph{Coding} 

In many research projects in nonlinear algebra the symbolic and the numerical approach are both needed for solving different parts of the same problem, and turn out to be highly interlinked. As a result, scientists frequently use different software and programming languages even for a single research project. An example of this is shown on the \MR page \texttt{\MYhref[burgundy]{https://mathrepo.mis.mpg.de/Landau/index.html}{Landau Discriminants}}, providing auxiliary material to \cite{mizera2021landau}.
This work applies methods from nonlinear algebra to the theory of scattering amplitudes in particle physics. In particular, these seemingly far-apart fields are practically and theoretically connected in the study of the so-called Landau equations. These are a set of polynomial equations determining allowed positions of singularities of Feynman integrals which arise in quantum field theory. The Landau discriminant parametrizes points for which the Landau equations have solutions. For each Feynman integral, the authors introduce the
Landau discriminant as a projective variety whose points are potential singularities of the integral. 
Nonlinear algebraists can then foster the understanding of this physical setting by studying geometric properties of the discriminant, such as irreducibility, dimension, and degree.

Together with the theoretical findings, they provide an implementation of their work. Their \texttt{Julia} package \texttt{Landau.jl} numerically computes defining equations of the Landau discriminant for some Feynman integrals which were previously out of reach. A tutorial illustrating how to use the package is provided as a \texttt{Jupyter} notebook, linked to the \MR page.
Further symbolic elimination methods are implemented in \texttt{Macaulay2} for computing Landau discriminants and are also presented and illustrated on the project page.
These approaches complement and enrich each other: the symbolic method provides reliable outputs but cannot deal with examples involving a high number of variables while the numerical computations are not exact but can be used to compute the desired equations in an efficient way. 

This project illustrates one big advantage of using \MR over different, say software-specific solutions: there is no constraint on choosing a particular programming language to work with, and the mathematician is free to do whatever their problem demands. Indeed, code fragments written in different languages can be easily combined and illustrated on the same webpage, together with information about the respective software version and hardware setup provided at the bottom of the project page. An additional written explanation between these code snippets greatly facilitates interoperability between different systems and improves reusability for the reader.

\begin{myfigure}
 \centering
 \includegraphics[width=\linewidth]{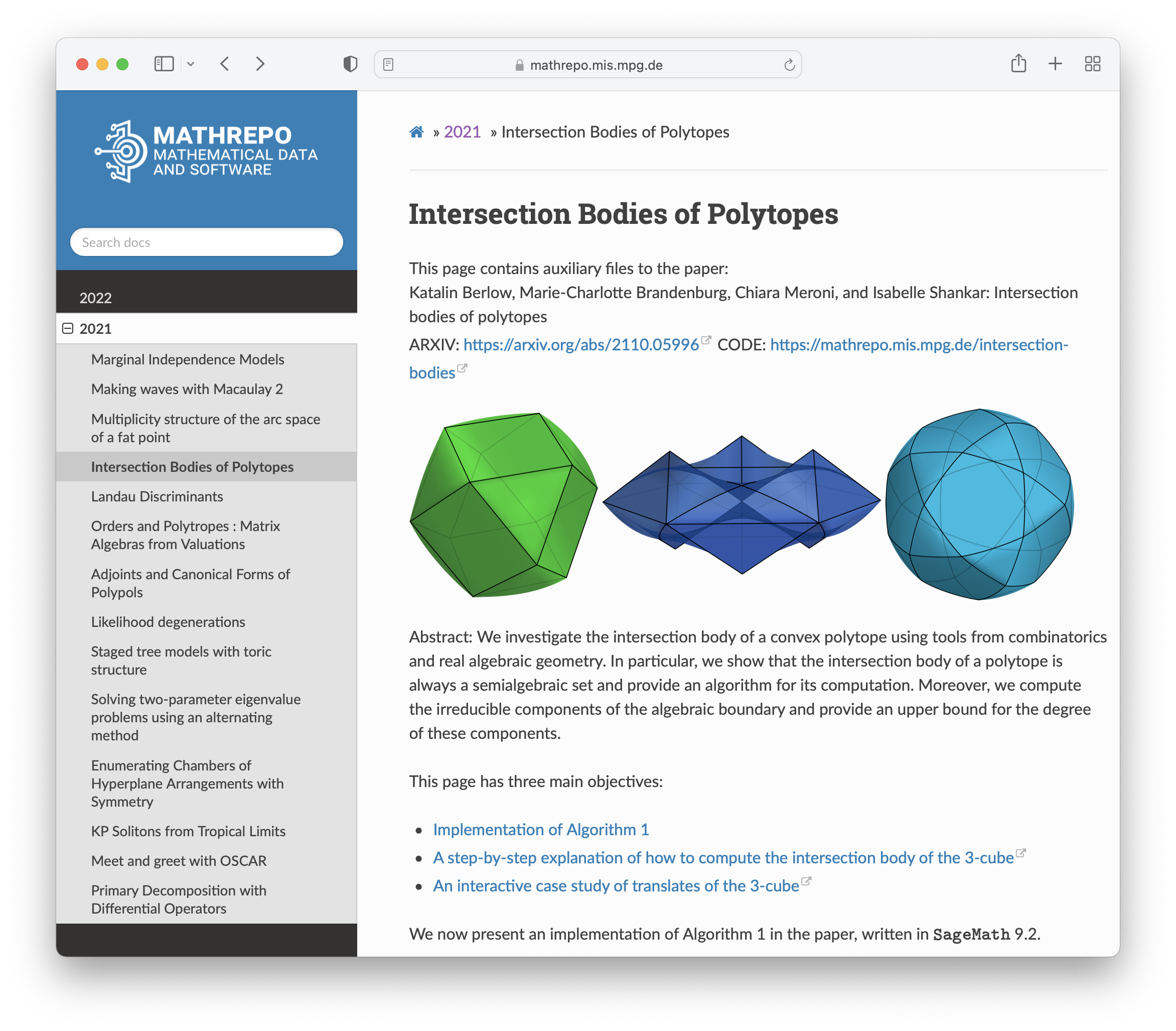}
 \captionof{figure}{A screenshot of the \MR subpage \texttt{Intersection Bodies of Polytopes}. 
 The \MR logo is displayed in the top-left corner, below that is the list of all contributions sorted by year.%
 }
 \label{figure}
\end{myfigure}

\paragraph{Visualization}
Another key benefit of putting research-data on \MR is that important features of the various programming languages, such as interactive visualizations, can be presented. The page \texttt{\MYhref[burgundy]{https://mathrepo.mis.mpg.de/intersection-bodies/index.html}{Intersection Bodies of Polytopes}}, illustrated in Fig.~\ref{figure}, makes use of this. In the corresponding article \cite{berlow2022intersection}, the authors investigate the intersection body of a convex polytope combining tools from combinatorics and from real algebraic geometry. In particular, they implement an algorithm for computing its algebraic boundary in \texttt{Sage} and \texttt{OSCAR}. Visitors of this page can both download the \texttt{Jupyter} notebooks to run the computations independently and also directly look at the code together with an interactive gallery of the output.

Finally, a step-by-step explanation of how to compute the intersection body of the three-dimensional cube, and an interactive case study of translates of the 3-cube are provided. In contrast to a paper publication, the \MR storage solution gives access to three-dimensional visualizations of these objects. In this way curvature, convexity, and all pieces of the algebraic boundary can be investigated.

\paragraph{Teaching}
A third type of contribution to \MR are presentations of problems solved in workshops or lectures. Indeed, the repository is also a possible useful tool for teaching purposes or for learning how to code in the languages previously specified. For instance, \texttt{\MYhref[burgundy]{https://mathrepo.mis.mpg.de/InvitationToNonlinearAlgebra/index.html}{Invitation to Nonlinear Algebra}} is a page where exercises, and examples dealing with polynomial rings, primary decomposition of ideals, tropical geometry, and tensors
presented in the textbook \cite{michalek2021invitation} are solved and explained using \texttt{Macaulay2} and \texttt{Polymake}. \medskip

The webpages from \MR we illustrate here are just some examples of many very different contributions included in the repository and they emphasize some of its main features. An autonomous navigation of the webpage is encouraged to get a more accurate idea of its possible types of use and to discover more interesting code presentations that are the building blocks of intuitions, examples, theorems, and their proofs.

\section{A discussion of FAIRness}

Nowadays there are general, not mathematics-specific guidelines for producing FAIR research data \cite[e.g.]{FAIRimplementation} and some of MaRDI's partnering research-data consortia in other disciplines are advanced in the process of implementing an apt infrastructure across their communities. Centrally, these implementations address \emph{automated} search and usage of research data. Especially with respect to findability and accessability, extensive metadata in some agreed-upon standard format is then essential. To implement interoperability,  it need be possible to automatically combine different digital resources, and for reusability automated decisions of relevance and legal terms need to be made. These developments were not known to the \MR group until the kickoff of MaRDI in 2021 and are, at the time of writing, not yet established in the mathematics community. However, all four principles have throughout been implicitly addressed from a human-user standpoint as follows.

\paragraph{Findability} Every subpage of \MR has a unique and persistent URL of the format \url{https://mathrepo.mis.mpg.de/<project-name>} assigned by its respective authors. Project names are encouraged to be telling and related to the title of the corresponding publication, if existent. The MPI MiS's internal library service provides up-to-date references for the latter. Vice versa, the link to the \MR project page can be explicitly stated in a paper publication.
The additional index on the repository's website which sorts entries by year and by programming language enables a direct as well as an associative findability on the webpage itself. These measures facilitate findability for humans.
The MPI MiS's front webpage contains a link to the repository and it is known to the most widely-used search engines, thus also improving machine findability. Persistence of the repository and its links is guaranteed by the MPI MiS's directors' decision and is independent of the respective project and author location. 

\paragraph{Accessability} All content available on \MR, both research data and their metadata, is freely accessible for all users from any location.
The current project-page template features the following recommendations: provide author names, a citation of a relevant publication, possibly \texttt{arXiv} link, an abstract, and the system setup (programming language, version, hardware) for computations, as well as a corresponding author for the project page itself. These minimum requirements establish a local metadata standard while allowing for a lot of flexibility in actual content.
The website is build on the \texttt{http} protocol, amenable to automated machine-search. The underlying \texttt{GitLab} is an open-source software. It provides an easy tracking tool, making historical changes and versions available for maintainers and contributors. These implementations make the website reasonably accessible without introducing a too technical setup.

\paragraph{Interoperability} 
On a theoretical level, \MR facilitates interoperability of its content by providing an abstract for each page. This embeds the employed mathematical terms into context and thus allows readers to translate the content into their own mathematical language, making it interoperable for humans. On the computational side, the more recent pages on \MR provide information on how to reproduce their respective system setup, facilitating interoperability on a technical level.

\paragraph{Reusability} Two implementations help readers to reproduce and reuse the \MR content. First, the abstract of each project page is often taken from a corresponding paper publication. It gives an indication as to which area of mathematics the research data belongs to and where to find additional background literature. Second, the pages created after 2018 each have a corresponding author named who can be contacted for direct information.
\medskip

The combination of these measures makes the data present on \MR FAIRer than research material which is simply put on personal homepages. \MR retains the low entry barrier and the flexibility of such user-implemented solutions, as shown in examples in the previous section, while providing long-term storage and independence of the authors' current academic location. To the best of our knowledge, this is a standalone feature of the repository in the mathematics community. However, many of these solutions are practical in a local sense and are not yet embedded in a national or even global infrastructure for research data.

\section{Current challenges}

At the time of writing, the \MR project standards do not follow a recommended protocol for mathematical research-data presentation and metadata supply. They rather reflect the maintainers' experience with documentation needed for the usage and maintenance of computer-algebra software. This approach has initially allowed for a swift setup of the repository. However, it does come with certain limitations, many of which became visible only with the growth of the repository. Now, challenges are threefold: inconsistent layout across the individual project pages, out-of-date or broken content, and large variations in mathematical quality.

The reasons and implications of these issues are numerous and the \MR contributions of the past five years clearly show a process of how priorities changed over time. For instance, the template for authors has improved with every handover between different sets of maintainers, providing increasingly refined metadata and new guidelines for presentation. As a result however, the depth of detail of the individual contributions' metadata and the overall layout of the \MR project pages is not uniform. This is now a hurdle to both findability and automated accessability, as well as a hurdle to reusability. For instance, metadata provided in different places across different pages encumbers findability for the human reader. Broken links to \eg \texttt{binder} notebooks are an obstacle to accessability and convey a general impression of content being out of date. Not sufficiently detailed references to employed programming languages hinder reproducibility on systems different from the authors'. Missing corresponding-author names in the early contributions impede reusability.

Even though training sessions at MPI MiS allow newcomers to learn to operate the \texttt{GitLab} setup, MathRepo's recommended project-page template is then often not followed in detail. This issue together with the great flexibility---to use \MR as a pure storage solution or, on the other extreme, as a platform for presenting teaching material---entail that contributions vary largely in quality.
But then a combination of large variations both in quality and in presentation can confuse readers and, worse, make them wary to trust the content. This is challenging from a user perspective.

From a maintainers' perspective, non-compliance with locally established standards in the new contributions and an increasing number of breakage in the past contributions hugely increase the workload that is needed to ensure some sort of consistency between the individual subpages of the repository.
Additionally, the user community has recently grown in numbers and has spread in mathematical diversity, making it hard for the maintainers to understand and judge new content. For scientists for whom \MR is only one of many projects and time commitment is limited, this is manageable only while the repository is still reasonably small. This issue actually is part of a bigger problem in a publishing culture which values code and software less than research data that comes in the form of paper publications.  But based on the past experiences and based on current research projects in the different working groups at MPI MiS, we predict more and more non-text research data and future growth of MathRepo. This snapshot marks a point at which it is still possible to tend to present issues before they become too large to manage.

\section{Outlook}

With increasing awareness of the importance of FAIR principles for sustainable research and with MaRDI launched, it is timely for us to address the above points of criticism.
In particular, we plan to tackle three main issues: trustworthiness of the repository, barriers to reusability, and forward-compatibility of the standards for embedding into coming MaRDI infrastructure.

With respect to trustworthiness, we see two possible ways for improvement of the status quo: to implement a user check for executability of code and for mathematical correctness of the research data similar to a small-scale peer-review system,
or to follow an established standard for research-data repositories and to apply for third-party certification, \eg with CoreTrustSeal\footnote{\url{https://www.coretrustseal.org}}.

This latter idea immediately leads to our point on reusability. Many seals of quality require a clear statement of the terms of use of a repository. However, at the time of writing this is not present in \MR and the lack of a license statement does not imply that data is automatically open access \cite[cf.]{FAIRimplementation}. 
For future FAIRness of the repository, it is thus mandatory to choose an appropriate license. This can either be by maintainer's choice or each user could be required to choose their own from a list of standard open-access licenses.

We aim to address the above issues and also facilitate a future embedding into MaRDI infrastructure by setting up a new template with metadata standards. These shall beforehand be discussed within the current MaRDI community. Compliance with the new standards shall then be achieved by embedding a section on research-data management into coming MathRepo-training sessions, stressing the importance of making your research FAIR.

\vfill\null
\columnbreak
\medskip
\bibliographystyle{alpha}
\bibliography{references}

\newcommand{\etalchar}[1]{$^{#1}$}
\begin{thebibliography}{JdMAJ{\etalchar{+}}20}

\bibitem[BBMS22]{berlow2022intersection}
Katalin Berlow, Marie-Charlotte Brandenburg, Chiara Meroni, and Isabelle
  Shankar.
\newblock Intersection bodies of polytopes.
\newblock {\em Beitr{\"a}ge zur Algebra und Geometrie/Contributions to Algebra
  and Geometry}, pages 1--21, 2022.

\bibitem[B{\c{C}}D{\etalchar{+}}21]{breiding2021nonlinear}
Paul Breiding, T{\"u}rk{\"u}~{\"O}zl{\"u}m {\c{C}}elik, Timothy Duff, Alexander
  Heaton, Aida Maraj, Anna-Laura Sattelberger, Lorenzo Venturello, and
  O{\u{g}}uzhan Y{\"u}r{\"u}k.
\newblock Nonlinear algebra and applications.
\newblock {\em Preprint available at \texttt{\emph{arXiv:2103.16300}}}, 2021.

\bibitem[BEO02]{smallgroups}
Hans~Ulrich Besche, Bettina Eick, and E.~A. O'Brien.
\newblock A millennium project: constructing small groups.
\newblock {\em Internat. J. Algebra Comput.}, 12(5):623--644, 2002.

\bibitem[GS21]{DMVM}
Christiane G\"orgen and Rainer Sinn.
\newblock Mathematik in der {Na}tionalen {F}orschungsdateninfrastruktur.
\newblock {\em Mitteilungen der Deutschen Mathematiker-Vereinigung},
  29(3):122--123, 2021.

\bibitem[HWS21]{NFDI}
Nathalie Hartl, Elena W{\"{o}}ssner, and York Sure{-}Vetter.
\newblock Nationale {F}orschungsdateninfrastruktur {(NFDI)}.
\newblock {\em Inform. Spektrum}, 44(5):370--373, 2021.

\bibitem[JdMAJ{\etalchar{+}}20]{FAIRimplementation}
Annika Jacobsen, Ricardo de~Miranda~Azevedo, Nick Juty, Dominique Batista,
  Simon Coles, Ronald Cornet, M{\'e}lanie Courtot, Merc{\`e} Crosas, Michel
  Dumontier, et~al.
\newblock {FAIR} principles: Interpretations and implementation considerations.
\newblock {\em Data Intelligence}, 2(1-2):10--29, 2020.

\bibitem[MS21]{michalek2021invitation}
Mateusz Micha{\l}ek and Bernd Sturmfels.
\newblock {\em Invitation to nonlinear algebra}, volume 211.
\newblock American Mathematical Soc., 2021.

\bibitem[MT21]{mizera2021landau}
Sebastian Mizera and Simon Telen.
\newblock Landau discriminants.
\newblock {\em {P}reprint available at \texttt{\emph{arXiv:2109.08036}}}, 2021.

\bibitem[Stu22]{sturmfels2021beyond}
Bernd Sturmfels.
\newblock Beyond linear algebra.
\newblock In {\em Proceedings of the International Congress of Mathematicians,
  St.\ Petersburg}, 2022.
\newblock Preprint available at \texttt{arXiv:2108.09494}.

\bibitem[WDA{\etalchar{+}}16]{FAIR}
Mark Wilkinson, Michel Dumontier, IJsbrand~Jan Aalbersberg, Gaby Appleton,
  et~al.
\newblock The {FAIR} guiding principles for scientific data management and
  stewardship.
\newblock {\em Scientific Data}, 3(160018), 2016.

\end{thebibliography}

\end{multicols}
\end{document}